\documentclass{article}
\usepackage{amsfonts}
\usepackage{amsmath}

\newtheorem{theorem}{Theorem}

\newtheorem{definition}[theorem]{Definition}

\input{tcilatex}
\begin{document}

\title{Intrinsic equations for a relaxed elastic line of second kind in
Minkowski 3-space}
\author{Ergin Bayram and Emin Kasap \\
Ondokuz May\i s University, Faculty of Arts and Sciences,\\
Mathematics Department, Samsun, Turkey}
\maketitle

\begin{abstract}
Let $\alpha $ be an arc on a connected oriented surface $S$ in Minkowski
3-space, parameterized by arc length $s$, with torsion $\tau $ and length $l$%
. The total square torsion $H$ of $\alpha $ is defined by $%
H=\int_{0}^{l}\tau ^{2}ds$. The arc is called a relaxed elastic line of
second kind if it is an extremal for the variational problem of minimizing
the value of $H$ within the family of all arcs of length $l$ on $S$ having
the same initial point and initial direction as $\alpha $. In this study, we
obtain differential equation and boundary conditions for a relaxed elastic
line of second kind on an oriented surface in Minkowski 3-space.
\end{abstract}

\section{Preliminaries and Introduction}

In this section, we give some fundamentals required for this paper.

\begin{definition}
$E^{n}\ $equipped with the metric%
\begin{equation*}
\left\langle u,w\right\rangle =-\sum_{i=1}^{\upsilon }u_{i}w_{i}+\sum_{j=\nu
+1}^{n}u_{j}w_{j,\ \ }\ \ \ u,w\in E^{n},\ 0\leq \nu \leq n,
\end{equation*}%
is called semi-Euclidean space and is denoted by $E_{\nu }^{n},$ where $\nu $
is called the index of the metric. For $n=3$, $E_{1}^{3}\ $is called
Minkowski 3-space \cite{oneill}.
\end{definition}

Let $\alpha \left( s\right) $ denote an arc on a connected oriented surface $%
S$ in $E_{1}^{3}$, parameterized by arc length $s$, $0\leq s\leq l$, with
curvature $\kappa \left( s\right) $ and torsion $\tau \left( s\right) $. Let
the energy density be given as some function of the curvature and torsion, $%
f\left( \kappa ,\ \tau \right) $. Then%
\begin{equation}
H=\int f\left( \kappa ,\ \tau \right) ds  \label{1}
\end{equation}%
is an Hamiltonian for curves \cite{hamilton}. Thus the following integral
can be taken as a special case of Hamiltonians for curves:%
\begin{equation}
H=\int \tau ^{2}ds.  \label{2.1}
\end{equation}%
Ekici and G\"{o}rg\"{u}l\"{u} \cite{ekici} handled the problem of
minimization of the integral $\int \kappa ^{2}\tau ds\ $in Minkowski
3-space. In \cite{ayhans} authors defined relaxed elastic line of second
kind on an oriented surface in Minkowski space and for a curve of this kind
lying on an oriented surface the Euler-Lagrange equations were derived. In
this paper, we give intrinsic equations for a curve of this kind.
Particularly, we obtain differential equation and boundary conditions for a
curve to be a relaxed elastic line of second kind on an oriented surface.

\begin{definition}
\bigskip The arc $\alpha $ is called a relaxed elastic line of second kind
in Minkowski 3-space if it is an extremal for the variational problem of
minimizing the value of $H$ within the family of all arcs of length $l$ on $%
S $ having the same initial point and initial direction as $\alpha $ in
Minkowski 3-space.
\end{definition}

In this study, we would like to calculate the intrinsic equations for the
curve $\alpha $ which is an extremal for $\left( \ref{1}\right) $. We shall
require that the coordinate functions of $S$ are smooth enough to have
partial derivatives and coordinate functions of $\alpha $, as functions of $%
s $, are smooth enough in these coordinates.

\begin{definition}
A tangent vector $v$ in $E_{1}^{3}$ is

spacelike if $\left\langle v,\ v\right\rangle $ $>0$ or $v=0$,

null if $\left\langle v,\ v\right\rangle =0$ and $v\neq 0$,

timelike if $\left\langle v,\ v\right\rangle $ $<0$.
\end{definition}

\begin{definition}
At a point $\alpha \left( s\right) $ of $\alpha $, let $T$ denote the unit
tangent vector to $\alpha $, $n$ the unit normal to $S$, and%
\begin{equation}
n\times T=\varepsilon Q\left( s\right) ,\ \ \varepsilon =\pm 1,  \label{3}
\end{equation}%
respectively. Then $\left\{ T,\ Q,\ n\right\} $ gives an orthonormal basis
in $E_{1}^{3}$. If $S$ is a spacelike surface then $T\times Q=n,\ Q\times
n=-T,\ n\times T=-Q.$ Similarly if $S$ is a\ timelike surface then $T\times
Q=-n,\ Q\times n=\pm T,\ n\times T=\mp Q$ $\cite{ugurlu}$.

\begin{theorem}
Let $S$ be a surface in $E_{1}^{3}$ and $\alpha \ $be a curve on $S$. The
analogue of the Frenet--Serret formulas is given by%
\begin{equation}
\left( 
\begin{array}{c}
T^{\prime } \\ 
Q^{\prime } \\ 
n^{\prime }%
\end{array}%
\right) =\left( 
\begin{array}{c}
\ \ \ \ \ 0\ \ \ \ \ \ \varepsilon _{2}\kappa _{g}\ \ \ \ \ \ \varepsilon
_{3}\kappa _{n} \\ 
-\varepsilon _{1}\kappa _{g}\ \ \ \ \ 0\ \ \ \ \ \ \ \ \ \ \varepsilon
_{3}\tau _{g} \\ 
-\varepsilon _{1}\kappa _{n}\ \ \ -\varepsilon _{2}\tau _{g}\ \ \ \ \ \ 0\ \
\ 
\end{array}%
\right) \left( 
\begin{array}{c}
T \\ 
Q \\ 
n%
\end{array}%
\right) ,  \label{2.2}
\end{equation}%
where $\varepsilon _{1}=\left\langle T,T\right\rangle $, $\varepsilon
_{2}=\left\langle Q,Q\right\rangle $,$\ \varepsilon _{3}=$ $\left\langle
n,n\right\rangle $. Here $k_{g}(s)=\left\langle T^{\prime
}(s),Q(s)\right\rangle $, $\tau _{g}(s)=\left\langle Q^{\prime
}(s),n(s)\right\rangle $ and $k_{n}(s)=\left\langle II\left( T\left(
s\right) ,T\left( s\right) \right) ,n\left( s\right) \right\rangle $%
\bigskip\ are geodesic curvature, geodesic torsion and normal curvature of $%
\alpha ,\ $respectively$.$
\end{theorem}

\begin{theorem}
Let $\alpha $ be any regular curve on a surface in $E_{1}^{3}.\ $Then we have%
\begin{equation*}
\kappa =\frac{\left\vert \left\vert \alpha ^{\prime }\times \alpha ^{\prime
\prime }\right\vert \right\vert }{\left\vert \left\vert \alpha ^{\prime
}\right\vert \right\vert ^{3}},\ \ \tau =\varepsilon _{1}\varepsilon _{2}%
\frac{\left\langle \alpha ^{\prime }\times \alpha ^{\prime \prime },\alpha
^{\prime \prime \prime }\right\rangle }{\left\vert \left\vert \alpha
^{\prime }\times \alpha ^{\prime \prime }\right\vert \right\vert ^{2}}.
\end{equation*}
\end{theorem}
\end{definition}

\section{Derivation of equations}

Suppose that $\alpha $ lies in a coordinate patch $\left( u,v\right)
\rightarrow x\left( u,v\right) \ $of a surface $S,\ $and let $x_{u}=\partial
x/\partial u,\ x_{v}=\partial x/\partial v.\ $Then $\alpha $ is expressed as

\begin{equation*}
\alpha \left( s\right) =x\left( u\left( s\right) ,v\left( s\right) \right)
,\ 0\leq s\leq l,
\end{equation*}%
with

\begin{equation*}
T\left( s\right) =\alpha ^{\prime }\left( s\right) =\frac{du}{ds}x_{u}+\frac{%
dv}{ds}x_{v}
\end{equation*}%
and

\begin{equation*}
Q\left( s\right) =p\left( s\right) x_{u}+q\left( s\right) x_{v}
\end{equation*}%
for suitable scalar functions $p\left( s\right) \ $and $q\left( s\right) .$

Now we will define variational fields for our problem. In order to obtain
variational arcs of length $l$, we need to extend $\alpha $ to an arc $%
\alpha ^{\ast }\left( s\right) $ defined for $0\leq s\leq l^{\ast },$ with $%
l^{\ast }\geq l$ but sufficiently close to $l$ so that $\alpha ^{\ast }$
lies in the coordinate patch. Let $\mu \left( s\right) ,\ 0\leq s\leq
l^{\ast },\ $be a scalar function of class $C^{2},\ $not vanishing
identically. Define%
\begin{equation*}
\eta \left( s\right) =\mu \left( s\right) p^{\ast }\left( s\right) ,\ \
\zeta \left( s\right) =\mu \left( s\right) q^{\ast }\left( s\right) .
\end{equation*}

Then

\begin{equation}
\eta \left( s\right) x_{u}+\zeta \left( s\right) x_{v}=\mu \left( s\right)
Q\left( s\right)  \label{2.4}
\end{equation}%
along $\alpha .\ $Also assume that

\begin{equation}
\mu \left( 0\right) =0,\ \mu ^{\prime }\left( 0\right) =0.  \label{2.5}
\end{equation}

Now define

\begin{equation}
\beta \left( \sigma ;t\right) =x\left( u\left( \sigma \right) +t\eta \left(
\sigma \right) ,v\left( \sigma \right) +t\zeta \left( \sigma \right) \right)
,  \label{2.6}
\end{equation}%
for $0\leq \sigma \leq $\bigskip $l^{\ast }.$ For $\left\vert t\right\vert
<\varepsilon _{1}\ $(where $\varepsilon _{1}>0\ $depends upon the choice of $%
\alpha ^{\ast }$ and of $\mu $), the point $\beta \left( \sigma ;t\right) $
lies in the coordinate patch. For fixed $t,\ \beta \left( \sigma ;t\right) \ 
$gives an arc with the same initial point and initial direction as $\alpha
,\ $because of $\left( \ref{2.5}\right) $. For $t=0,$ $\beta \left( \sigma
;0\right) \ $is the same as $\alpha ^{\ast }\ $and $\sigma \ $is arc length.
For $t\neq 0,\ $the parameter $\sigma \ $is not arc length in general.

For fixed $t,\ \left\vert t\right\vert <\varepsilon _{1},$ let $L^{\ast
}\left( t\right) $ denote the length of the arc $\beta \left( \sigma
;t\right) $, $0\leq \sigma \leq l^{\ast }.$ Then

\begin{equation}
L^{\ast }\left( t\right) =\int_{0}^{l}\sqrt{\left\vert \left\langle \frac{%
\partial \beta }{\partial \sigma },\frac{\partial \beta }{\partial \sigma }%
\right\rangle \right\vert }d\sigma  \label{2.7}
\end{equation}%
with

\begin{equation}
L^{\ast }\left( 0\right) =l^{\ast }>l.  \label{2.8}
\end{equation}

\bigskip By $\left( \ref{2.6}\right) $and $\left( \ref{2.7}\right) \ L^{\ast
}\left( t\right) \ $is continuous and differentiable in $t.\ $Particularly,
it follows from $\left( \ref{2.8}\right) $ that

\begin{equation}
L^{\ast }\left( t\right) >\frac{l+l^{\ast }}{2}>l\ \ \ \text{for \ }%
\left\vert t\right\vert <\varepsilon  \label{2.9}
\end{equation}%
for a suitable $\varepsilon \ $satisfying $0<\varepsilon \leq \varepsilon
_{1}.\ $Because of $\left( \ref{2.9}\right) \ $one can restrict $\beta
\left( \sigma ;t\right) \ ,\ 0\leq \left\vert t\right\vert <\varepsilon ,\ $%
to an arc of length $l$ by restricting the parameter $\sigma $ to an
interval $0\leq \sigma \leq \lambda \left( t\right) \leq l^{\ast }$ by
requiring

\begin{equation}
\int_{0}^{\lambda \left( t\right) }\sqrt{\left\vert \left\langle \frac{%
\partial \beta }{\partial \sigma },\frac{\partial \beta }{\partial \sigma }%
\right\rangle \right\vert }d\sigma =l.  \label{2.10}
\end{equation}

Note that $\lambda \left( 0\right) =l.$ The function $\lambda \left(
t\right) $ need not be determined explicitly, but we shall need

\begin{equation}
\left. \frac{d\lambda }{dt}\right\vert _{t=0}=\varepsilon
_{1}\int_{0}^{l}\mu \kappa _{g}ds.  \label{2.11}
\end{equation}

The proof of $\left( \ref{2.11}\right) $ and of other results will depend on
calculations from $\left( \ref{2.6}\right) $ such as

\begin{equation}
\left. \frac{\partial \beta }{\partial \sigma }\right\vert _{t=0}=T,\ \
0\leq s\leq l,  \label{2.12}
\end{equation}%
which gives

\begin{equation}
\left. \frac{\partial ^{2}\beta }{\partial \sigma ^{2}}\right\vert
_{t=0}=T^{\prime }=\varepsilon _{2}\kappa _{g}Q+\varepsilon _{3}\kappa _{n}n.
\label{2.13}
\end{equation}

Also

\begin{equation}
\left. \frac{\partial \beta }{\partial t}\right\vert _{t=0}=\mu Q
\label{2.14}
\end{equation}%
because of $\left( \ref{2.4}\right) $. Further differentiation of $\left( %
\ref{2.14}\right) $ gives

\begin{equation}
\left. \frac{\partial ^{2}\beta }{\partial t\partial \sigma }\right\vert
_{t=0}=\left. \frac{\partial ^{2}\beta }{\partial \sigma \partial t}%
\right\vert _{t=0}=\mu ^{\prime }Q+\mu Q^{\prime }=-\varepsilon _{1}\mu
\kappa _{g}T+\mu ^{\prime }Q+\varepsilon _{3}\mu \tau _{g}n  \label{2.15}
\end{equation}%
and using $\left( \ref{2.2}\right) $,%
\begin{eqnarray}
\left. \frac{\partial ^{3}\beta }{\partial t\partial \sigma ^{2}}\right\vert
_{t=0} &=&\left( -2\varepsilon _{1}\mu ^{\prime }\kappa _{g}-\varepsilon
_{1}\mu \kappa _{g}^{\prime }-\varepsilon _{1}\varepsilon _{3}\mu \kappa
_{n}\tau _{g}\right) T  \label{2.16} \\
&&+\left( \mu ^{\prime \prime }-\varepsilon _{1}\varepsilon _{2}\mu \kappa
_{g}^{2}-\varepsilon _{1}\varepsilon _{3}\mu \tau _{g}^{2}\right) Q  \notag
\\
&&+\left( 2\varepsilon _{3}\mu ^{\prime }\tau _{g}-\varepsilon
_{1}\varepsilon _{3}\mu \tau _{g}^{\prime }+\varepsilon _{1}\varepsilon
_{3}\mu \kappa _{g}\kappa _{n}\right) n.  \notag
\end{eqnarray}

Also using $\left( \ref{2.13}\right) $ we have%
\begin{eqnarray}
\left. \frac{\partial ^{3}\beta }{\partial \sigma ^{3}}\right\vert _{t=0}
&=&-\left( \varepsilon _{1}\varepsilon _{2}\kappa _{g}^{2}+\varepsilon
_{1}\varepsilon _{3}\kappa _{n}^{2}\right) T+\left( \varepsilon _{2}\kappa
_{g}^{\prime }-\varepsilon _{2}\varepsilon _{3}\kappa _{n}\tau _{g}\right) Q
\label{2.17} \\
&&+\left( \varepsilon _{3}\kappa _{n}^{\prime }+\varepsilon _{2}\varepsilon
_{3}\kappa _{g}\tau _{g}\right) n  \notag
\end{eqnarray}

and by $\left( \ref{2.14}\right) $

\begin{eqnarray}
\left. \frac{\partial ^{4}\beta }{\partial t\partial \sigma ^{3}}\right\vert
_{t=0} &=&\left( -3\varepsilon _{1}\mu ^{\prime \prime }\kappa
_{g}-3\varepsilon _{1}\mu ^{\prime }\kappa _{g}^{\prime }-\varepsilon
_{1}\mu \kappa _{g}^{\prime \prime }-\varepsilon _{1}\varepsilon _{3}\mu
\kappa _{n}^{\prime }\tau _{g}-2\varepsilon _{1}\varepsilon _{3}\mu \kappa
_{n}\tau _{g}^{\prime }\right.  \label{2.18} \\
&&\left. -3\varepsilon _{1}\varepsilon _{3}\mu ^{\prime }\kappa _{n}\tau
_{g}+\varepsilon _{3}\mu \kappa _{g}\kappa _{n}^{2}+\varepsilon _{2}\mu
\kappa _{g}^{3}+\varepsilon _{1}\varepsilon _{2}\varepsilon _{3}\mu \kappa
_{g}\tau _{g}^{2}\right) T  \notag \\
&&+\left( \mu ^{\prime \prime \prime }-3\varepsilon _{1}\varepsilon _{2}\mu
^{\prime }\kappa _{g}^{2}-3\varepsilon _{2}\varepsilon _{3}\mu ^{\prime
}\tau _{g}^{2}-3\varepsilon _{1}\varepsilon _{2}\mu \kappa _{g}\kappa
_{g}^{\prime }-3\varepsilon _{2}\varepsilon _{3}\mu \tau _{g}\tau
_{g}^{\prime }\right) Q  \notag \\
&&+\left( -3\varepsilon _{1}\varepsilon _{3}\mu ^{\prime }\kappa _{g}\kappa
_{n}-2\varepsilon _{1}\varepsilon _{3}\mu \kappa _{g}^{\prime }\kappa
_{n}-\varepsilon _{1}\varepsilon _{2}\varepsilon _{3}\mu \kappa _{g}^{2}\tau
_{g}+3\varepsilon _{3}\mu ^{\prime \prime }\tau _{g}\right.  \notag \\
&&\left. +3\varepsilon _{3}\mu ^{\prime }\tau _{g}^{\prime }+\varepsilon
_{3}\mu \tau _{g}^{\prime \prime }-\varepsilon _{1}\varepsilon _{3}\mu
\kappa _{g}\kappa _{n}^{\prime }-\varepsilon _{2}\mu \tau
_{g}^{3}-\varepsilon _{1}\mu \kappa _{n}^{2}\tau _{g}\right) n.  \notag
\end{eqnarray}

Now, let $H\left( t\right) $ denote the functional of a relaxed elastic line
of second kind for the arc $\beta \left( \sigma ;t\right) ,\ 0\leq \sigma
\leq \lambda \left( t\right) ,\ \left\vert t\right\vert <\varepsilon .$
Since, in general, $\sigma \ $is not the arc length for $t\neq 0$ functional 
$\left( \ref{2.1}\right) $ can be calculated as follows:%
\begin{equation*}
H\left( t\right) =\int_{0}^{\lambda \left( t\right) }\left( \frac{%
\left\langle \frac{\partial \beta }{\partial \sigma }\times \frac{\partial
^{2}\beta }{\partial \sigma ^{2}},\frac{\partial ^{3}\beta }{\partial \sigma
^{3}}\right\rangle }{\left\langle \frac{\partial \beta }{\partial \sigma },%
\frac{\partial \beta }{\partial \sigma }\right\rangle \left\langle \frac{%
\partial ^{2}\beta }{\partial \sigma ^{2}},\frac{\partial ^{2}\beta }{%
\partial \sigma ^{2}}\right\rangle -\left\langle \frac{\partial ^{2}\beta }{%
\partial \sigma ^{2}},\frac{\partial \beta }{\partial \sigma }\right\rangle
^{2}}\right) ^{2}d\sigma .
\end{equation*}

A necessary condition for $\alpha $ to be an extremal is that

\begin{equation*}
\left. \frac{dH}{dt}\right\vert _{t=0}=0
\end{equation*}%
for arbitrary $\mu $ satisfying $\left( \ref{2.5}\right) $. We have%
\begin{eqnarray*}
H^{\prime }\left( t\right) &=&\frac{d\lambda }{dt}\left\{ \left( \frac{%
\left\langle \frac{\partial \beta }{\partial \sigma }\times \frac{\partial
^{2}\beta }{\partial \sigma ^{2}},\frac{\partial ^{3}\beta }{\partial \sigma
^{3}}\right\rangle }{\left\langle \frac{\partial \beta }{\partial \sigma },%
\frac{\partial \beta }{\partial \sigma }\right\rangle \left\langle \frac{%
\partial ^{2}\beta }{\partial \sigma ^{2}},\frac{\partial ^{2}\beta }{%
\partial \sigma ^{2}}\right\rangle -\left\langle \frac{\partial ^{2}\beta }{%
\partial \sigma ^{2}},\frac{\partial \beta }{\partial \sigma }\right\rangle
^{2}}\right) ^{2}\right\} _{\sigma =\lambda \left( t\right) } \\
&&+2\int_{0}^{\lambda \left( t\right) }\frac{\partial }{\partial t}\left( 
\frac{\left\langle \frac{\partial \beta }{\partial \sigma }\times \frac{%
\partial ^{2}\beta }{\partial \sigma ^{2}},\frac{\partial ^{3}\beta }{%
\partial \sigma ^{3}}\right\rangle }{\left\langle \frac{\partial \beta }{%
\partial \sigma },\frac{\partial \beta }{\partial \sigma }\right\rangle
\left\langle \frac{\partial ^{2}\beta }{\partial \sigma ^{2}},\frac{\partial
^{2}\beta }{\partial \sigma ^{2}}\right\rangle -\left\langle \frac{\partial
^{2}\beta }{\partial \sigma ^{2}},\frac{\partial \beta }{\partial \sigma }%
\right\rangle ^{2}}\right) ^{2}d\sigma .
\end{eqnarray*}%
In calculating $dH/dt$; we give explicitly only terms that do not vanish for 
$t=0$. The omitted terms are those with factor

\begin{equation*}
\left\langle \frac{\partial ^{2}\beta }{\partial \sigma ^{2}},\frac{\partial
\beta }{\partial \sigma }\right\rangle
\end{equation*}%
which vanishes at $t=0$ because $\left\langle T,T^{\prime }\right\rangle =0.$
Thus, using $\left( \ref{2.10}-\ref{2.13}\right) $ and $\left( \ref{2.15}-%
\ref{2.18}\right) $ we get%
\begin{eqnarray*}
H^{\prime }\left( 0\right) &=&\varepsilon _{1}f^{2}\left( l\right)
\int_{0}^{l}\mu \kappa _{g}ds+2\int_{0}^{l}\frac{f\left( s\right) }{\left(
\varepsilon _{2}\kappa _{g}^{2}+\varepsilon _{3}\kappa _{n}^{2}\right) ^{2}}%
\left\{ \left( \varepsilon _{1}\varepsilon _{2}\kappa _{g}^{2}+\varepsilon
_{1}\varepsilon _{3}\kappa _{n}^{2}\right) \left[ \left( -3\kappa
_{g}^{2}\kappa _{n}^{\prime }\right. \right. \right. \\
&&-3\varepsilon _{2}\kappa _{g}^{3}\tau _{g}+3\kappa _{g}\kappa _{g}^{\prime
}\kappa _{n}-2\varepsilon _{3}\kappa _{g}\kappa _{n}^{2}\tau
_{g}+\varepsilon _{2}\kappa _{g}^{4}+\varepsilon _{3}\kappa _{g}^{2}\kappa
_{n}^{2}-\varepsilon _{2}\kappa _{n}^{\prime }\tau _{g}^{2}+\varepsilon
_{2}\varepsilon _{3}\kappa _{g}^{\prime }\tau _{g}^{\prime } \\
&&\left. +4\kappa _{n}\tau _{g}\tau _{g}^{\prime }-\varepsilon
_{2}\varepsilon _{3}\kappa _{g}\tau _{g}^{\prime \prime }+\varepsilon
_{1}\varepsilon _{2}\kappa _{g}\kappa _{n}^{2}\tau _{g}\right) +\mu ^{\prime
}\left( -4\kappa _{g}^{2}\kappa _{n}-\varepsilon _{1}\kappa
_{n}^{3}+2\varepsilon _{2}\varepsilon _{3}\kappa _{g}^{\prime }\tau
_{g}\right. \\
&&\left. +2\varepsilon _{2}\kappa _{n}\tau _{g}^{2}-3\varepsilon
_{2}\varepsilon _{3}\kappa _{g}\tau _{g}^{\prime }+3\kappa _{g}^{2}\kappa
_{n}+3\varepsilon _{2}\kappa _{n}\tau _{g}^{2}\right) +\mu ^{\prime \prime
}\left( -\varepsilon _{3}\kappa _{n}^{\prime }-4\varepsilon _{2}\varepsilon
_{3}\kappa _{g}\tau _{g}\right) \\
&&+\varepsilon _{3}\mu ^{\prime \prime \prime }\kappa _{n}-\left(
-\varepsilon _{2}\kappa _{g}\left( \varepsilon _{3}\kappa _{n}^{\prime
}+\varepsilon _{2}\varepsilon _{3}\kappa _{g}\tau _{g}\right) +\varepsilon
_{3}\kappa _{n}\left( \varepsilon _{2}\kappa _{g}^{\prime }-\varepsilon
_{2}\varepsilon _{3}\kappa _{n}\tau _{g}\right) \right) \left( \mu \left(
-4\varepsilon _{2}\kappa _{g}^{3}\right. \right. \\
&&\left. \left. \left. -4\varepsilon _{3}\kappa _{g}\kappa _{n}^{2}+2\kappa
_{g}\tau _{g}^{2}+2\varepsilon _{1}\varepsilon _{3}\kappa _{n}\tau
_{g}^{\prime }\right) +4\mu ^{\prime }\varepsilon _{1}\varepsilon _{3}\kappa
_{n}\tau _{g}+2\mu ^{\prime \prime }\varepsilon _{1}\kappa _{g}\right)
\right\} ds \\
&=&\int_{0}^{l}\mu \left( \varepsilon _{1}f^{2}\left( l\right) \kappa
_{g}\right) ds+\int_{0}^{l}\mu \left\{ \frac{2f\left( s\right) }{\left(
\varepsilon _{2}\kappa _{g}^{2}+\varepsilon _{3}\kappa _{n}^{2}\right) ^{2}}%
\left[ \left( \varepsilon _{1}\varepsilon _{2}\kappa _{g}^{2}+\varepsilon
_{1}\varepsilon _{3}\kappa _{n}^{2}\right) \left( -3\kappa _{g}^{2}\kappa
_{n}^{\prime }\right. \right. \right. \\
&&-3\varepsilon _{2}\kappa _{g}^{3}\tau _{g}+3\kappa _{g}\kappa _{g}^{\prime
}\kappa _{n}-2\varepsilon _{3}\kappa _{g}\kappa _{n}^{2}\tau
_{g}+\varepsilon _{2}\kappa _{g}^{4}+\varepsilon _{3}\kappa _{g}^{2}\kappa
_{n}^{2}-\varepsilon _{2}\kappa _{n}^{\prime }\tau _{g}^{2}+\varepsilon
_{2}\varepsilon _{3}\kappa _{g}^{\prime }\tau _{g}^{\prime } \\
&&\left. +4\kappa _{n}\tau _{g}\tau _{g}^{\prime }-\varepsilon
_{2}\varepsilon _{3}\kappa _{g}\tau _{g}^{\prime \prime }+\varepsilon
_{1}\varepsilon _{2}\kappa _{g}\kappa _{n}^{2}\tau _{g}\right) -\left(
-\varepsilon _{2}\kappa _{g}\left( \varepsilon _{3}\kappa _{n}^{\prime
}+\varepsilon _{2}\varepsilon _{3}\kappa _{g}\tau _{g}\right) \right. \\
&&\left. \left. \left. +\varepsilon _{3}\kappa _{n}\left( \varepsilon
_{2}\kappa _{g}^{\prime }-\varepsilon _{2}\varepsilon _{3}\kappa _{n}\tau
_{g}\right) \right) \left( -4\varepsilon _{2}\kappa _{g}^{3}-4\varepsilon
_{3}\kappa _{g}\kappa _{n}^{2}+2\kappa _{g}\tau _{g}^{2}+2\varepsilon
_{1}\varepsilon _{3}\kappa _{n}\tau _{g}^{\prime }\right) \right] \right\} ds
\\
&&+\int_{0}^{l}\mu ^{\prime }\left\{ \frac{2f\left( s\right) }{\left(
\varepsilon _{2}\kappa _{g}^{2}+\varepsilon _{3}\kappa _{n}^{2}\right) ^{2}}%
\left[ \left( \varepsilon _{1}\varepsilon _{2}\kappa _{g}^{2}+\varepsilon
_{1}\varepsilon _{3}\kappa _{n}^{2}\right) \left( -4\kappa _{g}^{2}\kappa
_{n}-\varepsilon _{1}\kappa _{n}^{3}+2\varepsilon _{2}\varepsilon _{3}\kappa
_{g}^{\prime }\tau _{g}\right. \right. \right. \\
&&\left. +2\varepsilon _{2}\kappa _{n}\tau _{g}^{2}-3\varepsilon
_{2}\varepsilon _{3}\kappa _{g}\tau _{g}^{\prime }+3\kappa _{g}^{2}\kappa
_{n}+3\varepsilon _{2}\kappa _{n}\tau _{g}^{2}\right) -4\varepsilon
_{1}\varepsilon _{3}\kappa _{n}\tau _{g}\left( -\varepsilon _{2}\kappa
_{g}\left( \varepsilon _{3}\kappa _{n}^{\prime }+\varepsilon _{2}\varepsilon
_{3}\kappa _{g}\tau _{g}\right) \right. \\
&&\left. \left. \left. +\varepsilon _{3}\kappa _{n}\left( \varepsilon
_{2}\kappa _{g}^{\prime }-\varepsilon _{2}\varepsilon _{3}\kappa _{n}\tau
_{g}\right) \right) \right] \right\} ds \\
&&+\int_{0}^{l}\mu ^{\prime \prime }\left\{ \frac{2f\left( s\right) }{\left(
\varepsilon _{2}\kappa _{g}^{2}+\varepsilon _{3}\kappa _{n}^{2}\right) ^{2}}%
\left[ \left( \varepsilon _{1}\varepsilon _{2}\kappa _{g}^{2}+\varepsilon
_{1}\varepsilon _{3}\kappa _{n}^{2}\right) \left( -\varepsilon _{3}\kappa
_{n}^{\prime }-4\varepsilon _{2}\varepsilon _{3}\kappa _{g}\tau _{g}\right)
\right. \right. \\
&&\left. \left. -2\varepsilon _{1}\kappa _{g}\left( -\varepsilon _{2}\kappa
_{g}\left( \varepsilon _{3}\kappa _{n}^{\prime }+\varepsilon _{2}\varepsilon
_{3}\kappa _{g}\tau _{g}\right) +\varepsilon _{3}\kappa _{n}\left(
\varepsilon _{2}\kappa _{g}^{\prime }-\varepsilon _{2}\varepsilon _{3}\kappa
_{n}\tau _{g}\right) \right) \right] \right\} ds \\
&&+\int_{0}^{l}\mu ^{\prime \prime \prime }\left\{ \frac{2\varepsilon
_{3}\kappa _{n}f\left( s\right) }{\left( \varepsilon _{2}\kappa
_{g}^{2}+\varepsilon _{3}\kappa _{n}^{2}\right) ^{2}}\right\} ds,
\end{eqnarray*}%
where%
\begin{equation*}
f\left( s\right) =\frac{-\varepsilon _{2}\kappa _{g}\left( \varepsilon
_{3}\kappa _{n}^{\prime }+\varepsilon _{2}\varepsilon _{3}\kappa _{g}\tau
_{g}\right) +\varepsilon _{3}\kappa _{n}\left( \varepsilon _{2}\kappa
_{g}^{\prime }-\varepsilon _{2}\varepsilon _{3}\kappa _{n}\tau _{g}\right) }{%
\varepsilon _{1}\left( \varepsilon _{2}\kappa _{g}^{2}+\varepsilon
_{3}\kappa _{n}^{2}\right) }.
\end{equation*}%
However, using integration by parts and $\left( \ref{2.5}\right) \ $we get%
\begin{eqnarray*}
&&\int_{0}^{l}\mu ^{\prime }\left\{ \frac{2f\left( s\right) }{\left(
\varepsilon _{2}\kappa _{g}^{2}+\varepsilon _{3}\kappa _{n}^{2}\right) ^{2}}%
\left[ \left( \varepsilon _{1}\varepsilon _{2}\kappa _{g}^{2}+\varepsilon
_{1}\varepsilon _{3}\kappa _{n}^{2}\right) \left( -4\kappa _{g}^{2}\kappa
_{n}-\varepsilon _{1}\kappa _{n}^{3}+2\varepsilon _{2}\varepsilon _{3}\kappa
_{g}^{\prime }\tau _{g}\right. \right. \right. \\
&&\left. +2\varepsilon _{2}\kappa _{n}\tau _{g}^{2}-3\varepsilon
_{2}\varepsilon _{3}\kappa _{g}\tau _{g}^{\prime }+3\kappa _{g}^{2}\kappa
_{n}+3\varepsilon _{2}\kappa _{n}\tau _{g}^{2}\right) -4\varepsilon
_{1}\varepsilon _{3}\kappa _{n}\tau _{g}\left( -\varepsilon _{2}\kappa
_{g}\left( \varepsilon _{3}\kappa _{n}^{\prime }+\varepsilon _{2}\varepsilon
_{3}\kappa _{g}\tau _{g}\right) \right. \\
&&\left. \left. \left. +\varepsilon _{3}\kappa _{n}\left( \varepsilon
_{2}\kappa _{g}^{\prime }-\varepsilon _{2}\varepsilon _{3}\kappa _{n}\tau
_{g}\right) \right) \right] \right\} ds=\mu \left( l\right) \left\{ \frac{%
2f\left( s\right) }{\left( \varepsilon _{2}\kappa _{g}^{2}+\varepsilon
_{3}\kappa _{n}^{2}\right) ^{2}}\left[ \left( \varepsilon _{1}\varepsilon
_{2}\kappa _{g}^{2}+\varepsilon _{1}\varepsilon _{3}\kappa _{n}^{2}\right)
\left( -4\kappa _{g}^{2}\kappa _{n}\right. \right. \right. \\
&&\left. -\varepsilon _{1}\kappa _{n}^{3}+2\varepsilon _{2}\varepsilon
_{3}\kappa _{g}^{\prime }\tau _{g}+2\varepsilon _{2}\kappa _{n}\tau
_{g}^{2}-3\varepsilon _{2}\varepsilon _{3}\kappa _{g}\tau _{g}^{\prime
}+3\kappa _{g}^{2}\kappa _{n}+3\varepsilon _{2}\kappa _{n}\tau
_{g}^{2}\right) \\
&&\left. \left. -4\varepsilon _{1}\varepsilon _{3}\kappa _{n}\tau _{g}\left(
-\varepsilon _{2}\kappa _{g}\left( \varepsilon _{3}\kappa _{n}^{\prime
}+\varepsilon _{2}\varepsilon _{3}\kappa _{g}\tau _{g}\right) +\varepsilon
_{3}\kappa _{n}\left( \varepsilon _{2}\kappa _{g}^{\prime }-\varepsilon
_{2}\varepsilon _{3}\kappa _{n}\tau _{g}\right) \right) \right] \right\}
_{s=l} \\
&&-\int_{0}^{l}\mu \left\{ \frac{2f\left( s\right) }{\left( \varepsilon
_{2}\kappa _{g}^{2}+\varepsilon _{3}\kappa _{n}^{2}\right) ^{2}}\left[
\left( \varepsilon _{1}\varepsilon _{2}\kappa _{g}^{2}+\varepsilon
_{1}\varepsilon _{3}\kappa _{n}^{2}\right) \left( -4\kappa _{g}^{2}\kappa
_{n}-\varepsilon _{1}\kappa _{n}^{3}+2\varepsilon _{2}\varepsilon _{3}\kappa
_{g}^{\prime }\tau _{g}\right. \right. \right. \\
&&\left. +2\varepsilon _{2}\kappa _{n}\tau _{g}^{2}-3\varepsilon
_{2}\varepsilon _{3}\kappa _{g}\tau _{g}^{\prime }+3\kappa _{g}^{2}\kappa
_{n}+3\varepsilon _{2}\kappa _{n}\tau _{g}^{2}\right) -4\varepsilon
_{1}\varepsilon _{3}\kappa _{n}\tau _{g}\left( -\varepsilon _{2}\kappa
_{g}\left( \varepsilon _{3}\kappa _{n}^{\prime }+\varepsilon _{2}\varepsilon
_{3}\kappa _{g}\tau _{g}\right) \right. \\
&&\left. \left. \left. +\varepsilon _{3}\kappa _{n}\left( \varepsilon
_{2}\kappa _{g}^{\prime }-\varepsilon _{2}\varepsilon _{3}\kappa _{n}\tau
_{g}\right) \right) \right] \right\} ^{\prime }ds
\end{eqnarray*}%
and%
\begin{eqnarray*}
&&\int_{0}^{l}\mu ^{\prime \prime }\left\{ \frac{2f\left( s\right) }{\left(
\varepsilon _{2}\kappa _{g}^{2}+\varepsilon _{3}\kappa _{n}^{2}\right) ^{2}}%
\left[ \left( \varepsilon _{1}\varepsilon _{2}\kappa _{g}^{2}+\varepsilon
_{1}\varepsilon _{3}\kappa _{n}^{2}\right) \left( -\varepsilon _{3}\kappa
_{n}^{\prime }-4\varepsilon _{2}\varepsilon _{3}\kappa _{g}\tau _{g}\right)
\right. \right. \\
&&\left. \left. -2\varepsilon _{1}\kappa _{g}\left( -\varepsilon _{2}\kappa
_{g}\left( \varepsilon _{3}\kappa _{n}^{\prime }+\varepsilon _{2}\varepsilon
_{3}\kappa _{g}\tau _{g}\right) +\varepsilon _{3}\kappa _{n}\left(
\varepsilon _{2}\kappa _{g}^{\prime }-\varepsilon _{2}\varepsilon _{3}\kappa
_{n}\tau _{g}\right) \right) \right] \right\} ds= \\
&&\mu ^{\prime }\left( l\right) \left\{ \frac{2f\left( s\right) }{\left(
\varepsilon _{2}\kappa _{g}^{2}+\varepsilon _{3}\kappa _{n}^{2}\right) ^{2}}%
\left[ \left( \varepsilon _{1}\varepsilon _{2}\kappa _{g}^{2}+\varepsilon
_{1}\varepsilon _{3}\kappa _{n}^{2}\right) \left( -\varepsilon _{3}\kappa
_{n}^{\prime }-4\varepsilon _{2}\varepsilon _{3}\kappa _{g}\tau _{g}\right)
\right. \right. \\
&&\left. \left. -2\varepsilon _{1}\kappa _{g}\left( -\varepsilon _{2}\kappa
_{g}\left( \varepsilon _{3}\kappa _{n}^{\prime }+\varepsilon _{2}\varepsilon
_{3}\kappa _{g}\tau _{g}\right) +\varepsilon _{3}\kappa _{n}\left(
\varepsilon _{2}\kappa _{g}^{\prime }-\varepsilon _{2}\varepsilon _{3}\kappa
_{n}\tau _{g}\right) \right) \right] \right\} _{s=l} \\
&&-\mu \left( l\right) \left\{ \frac{2f\left( s\right) }{\left( \varepsilon
_{2}\kappa _{g}^{2}+\varepsilon _{3}\kappa _{n}^{2}\right) ^{2}}\left[
\left( \varepsilon _{1}\varepsilon _{2}\kappa _{g}^{2}+\varepsilon
_{1}\varepsilon _{3}\kappa _{n}^{2}\right) \left( -\varepsilon _{3}\kappa
_{n}^{\prime }-4\varepsilon _{2}\varepsilon _{3}\kappa _{g}\tau _{g}\right)
\right. \right. \\
&&\left. \left. -2\varepsilon _{1}\kappa _{g}\left( -\varepsilon _{2}\kappa
_{g}\left( \varepsilon _{3}\kappa _{n}^{\prime }+\varepsilon _{2}\varepsilon
_{3}\kappa _{g}\tau _{g}\right) +\varepsilon _{3}\kappa _{n}\left(
\varepsilon _{2}\kappa _{g}^{\prime }-\varepsilon _{2}\varepsilon _{3}\kappa
_{n}\tau _{g}\right) \right) \right] \right\} _{s=l}^{\prime } \\
&&+\int_{0}^{l}\mu \left\{ \frac{2f\left( s\right) }{\left( \varepsilon
_{2}\kappa _{g}^{2}+\varepsilon _{3}\kappa _{n}^{2}\right) ^{2}}\left[
\left( \varepsilon _{1}\varepsilon _{2}\kappa _{g}^{2}+\varepsilon
_{1}\varepsilon _{3}\kappa _{n}^{2}\right) \left( -\varepsilon _{3}\kappa
_{n}^{\prime }-4\varepsilon _{2}\varepsilon _{3}\kappa _{g}\tau _{g}\right)
\right. \right. \\
&&\left. \left. -2\varepsilon _{1}\kappa _{g}\left( -\varepsilon _{2}\kappa
_{g}\left( \varepsilon _{3}\kappa _{n}^{\prime }+\varepsilon _{2}\varepsilon
_{3}\kappa _{g}\tau _{g}\right) +\varepsilon _{3}\kappa _{n}\left(
\varepsilon _{2}\kappa _{g}^{\prime }-\varepsilon _{2}\varepsilon _{3}\kappa
_{n}\tau _{g}\right) \right) \right] \right\} ^{\prime \prime }ds
\end{eqnarray*}%
and%
\begin{eqnarray*}
\int_{0}^{l}\mu ^{\prime \prime \prime }\left\{ \frac{2\varepsilon
_{3}\kappa _{n}f\left( s\right) }{\left( \varepsilon _{2}\kappa
_{g}^{2}+\varepsilon _{3}\kappa _{n}^{2}\right) ^{2}}\right\} ds &=&\mu
^{\prime \prime }\left( l\right) \left\{ \frac{2\varepsilon _{3}\kappa
_{n}f\left( s\right) }{\left( \varepsilon _{2}\kappa _{g}^{2}+\varepsilon
_{3}\kappa _{n}^{2}\right) ^{2}}\right\} _{s=l}-\mu ^{\prime \prime }\left(
0\right) \left\{ \frac{2\varepsilon _{3}\kappa _{n}f\left( s\right) }{\left(
\varepsilon _{2}\kappa _{g}^{2}+\varepsilon _{3}\kappa _{n}^{2}\right) ^{2}}%
\right\} _{s=0} \\
&&+\mu ^{\prime }\left( l\right) \left\{ \frac{2\varepsilon _{3}\kappa
_{n}f\left( s\right) }{\left( \varepsilon _{2}\kappa _{g}^{2}+\varepsilon
_{3}\kappa _{n}^{2}\right) ^{2}}\right\} _{s=l}^{\prime }+\mu \left(
l\right) \left\{ \frac{2\varepsilon _{3}\kappa _{n}f\left( s\right) }{\left(
\varepsilon _{2}\kappa _{g}^{2}+\varepsilon _{3}\kappa _{n}^{2}\right) ^{2}}%
\right\} _{s=l}^{\prime \prime } \\
&&-\int_{0}^{l}\mu \left\{ \frac{2\varepsilon _{3}\kappa _{n}f\left(
s\right) }{\left( \varepsilon _{2}\kappa _{g}^{2}+\varepsilon _{3}\kappa
_{n}^{2}\right) ^{2}}\right\} ^{\prime \prime \prime }ds.
\end{eqnarray*}%
Thus $H^{\prime }\left( 0\right) \ $can be written as%
\begin{eqnarray*}
H^{\prime }\left( 0\right) &=&\int_{0}^{l}\mu \left\{ \varepsilon _{1}\kappa
_{g}f^{2}\left( l\right) +\frac{2f\left( s\right) }{\left( \varepsilon
_{2}\kappa _{g}^{2}+\varepsilon _{3}\kappa _{n}^{2}\right) ^{2}}\left(
\left( \varepsilon _{1}\varepsilon _{2}\kappa _{g}^{2}+\varepsilon
_{1}\varepsilon _{3}\kappa _{n}^{2}\right) \left[ \left( -3\kappa
_{g}^{2}\kappa _{n}^{\prime }\right. \right. \right. \right. \\
&&-3\varepsilon _{2}\kappa _{g}^{3}\tau _{g}+3\kappa _{g}\kappa _{g}^{\prime
}\kappa _{n}-2\varepsilon _{3}\kappa _{g}\kappa _{n}^{2}\tau
_{g}+\varepsilon _{2}\kappa _{g}^{4}+\varepsilon _{3}\kappa _{g}^{2}\kappa
_{n}^{2}-\varepsilon _{2}\kappa _{n}^{\prime }\tau _{g}^{2}+\varepsilon
_{2}\varepsilon _{3}\kappa _{g}^{\prime }\tau _{g}^{\prime } \\
&&\left. +4\kappa _{n}\tau _{g}\tau _{g}^{\prime }-\varepsilon
_{2}\varepsilon _{3}\kappa _{g}\tau _{g}^{\prime \prime }+\varepsilon
_{1}\varepsilon _{2}\kappa _{g}\kappa _{n}^{2}\tau _{g}\right) -\left(
-\varepsilon _{2}\kappa _{g}\left( \varepsilon _{3}\kappa _{n}^{\prime
}+\varepsilon _{2}\varepsilon _{3}\kappa _{g}\tau _{g}\right) \right. \\
&&\left. \left. \left. +\varepsilon _{3}\kappa _{n}\left( \varepsilon
_{2}\kappa _{g}^{\prime }-\varepsilon _{2}\varepsilon _{3}\kappa _{n}\tau
_{g}\right) \right) \left( -4\varepsilon _{2}\kappa _{g}^{3}-4\varepsilon
_{3}\kappa _{g}\kappa _{n}^{2}+2\kappa _{g}\tau _{g}^{2}+2\varepsilon
_{1}\varepsilon _{3}\kappa _{n}\tau _{g}^{\prime }\right) \right] \right) \\
&&-\left[ \frac{2f\left( s\right) }{\left( \varepsilon _{2}\kappa
_{g}^{2}+\varepsilon _{3}\kappa _{n}^{2}\right) ^{2}}\left[ \left(
\varepsilon _{1}\varepsilon _{2}\kappa _{g}^{2}+\varepsilon _{1}\varepsilon
_{3}\kappa _{n}^{2}\right) \left( -4\kappa _{g}^{2}\kappa _{n}-\varepsilon
_{1}\kappa _{n}^{3}+2\varepsilon _{2}\varepsilon _{3}\kappa _{g}^{\prime
}\tau _{g}\right. \right. \right. \\
&&\left. +2\varepsilon _{2}\kappa _{n}\tau _{g}^{2}-3\varepsilon
_{2}\varepsilon _{3}\kappa _{g}\tau _{g}^{\prime }+3\kappa _{g}^{2}\kappa
_{n}+3\varepsilon _{2}\kappa _{n}\tau _{g}^{2}\right) -4\varepsilon
_{1}\varepsilon _{2}\kappa _{n}\tau _{g}\left( -\varepsilon _{2}\kappa
_{g}\left( \varepsilon _{3}\kappa _{n}^{\prime }+\varepsilon _{2}\varepsilon
_{3}\kappa _{g}\tau _{g}\right) \right. \\
&&\left. \left. \left. +\varepsilon _{3}\kappa _{n}\left( \varepsilon
_{2}\kappa _{g}^{\prime }-\varepsilon _{2}\varepsilon _{3}\kappa _{n}\tau
_{g}\right) \right) \right] \right] ^{\prime }+\left[ \frac{2f\left(
s\right) }{\left( \varepsilon _{2}\kappa _{g}^{2}+\varepsilon _{3}\kappa
_{n}^{2}\right) ^{2}}\left[ \left( \varepsilon _{1}\varepsilon _{2}\kappa
_{g}^{2}+\varepsilon _{1}\varepsilon _{3}\kappa _{n}^{2}\right) \left(
-\varepsilon _{3}\kappa _{n}^{\prime }\right. \right. \right. \\
&&\left. \left. \left. -4\varepsilon _{2}\varepsilon _{3}\kappa _{g}\tau
_{g}\right) -2\varepsilon _{1}\kappa _{g}\left( -\varepsilon _{2}\kappa
_{g}\left( \varepsilon _{3}\kappa _{n}^{\prime }+\varepsilon _{2}\varepsilon
_{3}\kappa _{g}\tau _{g}\right) +\varepsilon _{3}\kappa _{n}\left(
\varepsilon _{2}\kappa _{g}^{\prime }-\varepsilon _{2}\varepsilon _{3}\kappa
_{n}\tau _{g}\right) \right) \right] \right] ^{\prime \prime } \\
&&\left. -\left[ \frac{2\varepsilon _{3}\kappa _{n}f\left( s\right) }{\left(
\varepsilon _{2}\kappa _{g}^{2}+\varepsilon _{3}\kappa _{n}^{2}\right) ^{2}}%
\right] ^{\prime \prime \prime }\right\} ds \\
&&+\mu \left( l\right) \left\{ \frac{2f\left( l\right) }{\left( \varepsilon
_{2}\kappa _{g}^{2}\left( l\right) +\varepsilon _{3}\kappa _{n}^{2}\left(
l\right) \right) ^{2}}\left[ \left( \varepsilon _{1}\varepsilon _{2}\kappa
_{g}^{2}+\varepsilon _{1}\varepsilon _{3}\kappa _{n}^{2}\right) \left(
-4\kappa _{g}^{2}\kappa _{n}-\varepsilon _{1}\kappa _{n}^{3}+2\varepsilon
_{2}\varepsilon _{3}\kappa _{g}^{\prime }\tau _{g}\right. \right. \right. \\
&&\left. +2\varepsilon _{2}\kappa _{n}\tau _{g}^{2}-3\varepsilon
_{2}\varepsilon _{3}\kappa _{g}\tau _{g}^{\prime }+3\kappa _{g}^{2}\kappa
_{n}+3\varepsilon _{2}\kappa _{n}\tau _{g}^{2}\right) -4\varepsilon
_{1}\varepsilon _{2}\kappa _{n}\tau _{g}\left( -\varepsilon _{2}\kappa
_{g}\left( \varepsilon _{3}\kappa _{n}^{\prime }+\varepsilon _{2}\varepsilon
_{3}\kappa _{g}\tau _{g}\right) \right. \\
&&\left. \left. +\varepsilon _{3}\kappa _{n}\left( \varepsilon _{2}\kappa
_{g}^{\prime }-\varepsilon _{2}\varepsilon _{3}\kappa _{n}\tau _{g}\right)
\right) \right] _{s=l}-\left( \frac{2f\left( s\right) }{\left( \varepsilon
_{2}\kappa _{g}^{2}+\varepsilon _{3}\kappa _{n}^{2}\right) ^{2}}\left[
\left( \varepsilon _{1}\varepsilon _{2}\kappa _{g}^{2}+\varepsilon
_{1}\varepsilon _{3}\kappa _{n}^{2}\right) \right. \right. \\
&&x\left( -\varepsilon _{3}\kappa _{n}^{\prime }-4\varepsilon
_{2}\varepsilon _{3}\kappa _{g}\tau _{g}\right) -2\varepsilon _{1}\kappa
_{g}\left( -\varepsilon _{2}\kappa _{g}\left( \varepsilon _{3}\kappa
_{n}^{\prime }+\varepsilon _{2}\varepsilon _{3}\kappa _{g}\tau _{g}\right)
\right. \\
&&\left. \left. \left. \left. +\varepsilon _{3}\kappa _{n}\left( \varepsilon
_{2}\kappa _{g}^{\prime }-\varepsilon _{2}\varepsilon _{3}\kappa _{n}\tau
_{g}\right) \right) \right] \right) _{s=l}^{\prime }+\left( \frac{%
2\varepsilon _{3}\kappa _{n}f\left( s\right) }{\left( \varepsilon _{2}\kappa
_{g}^{2}+\varepsilon _{3}\kappa _{n}^{2}\right) ^{2}}\right) _{s=l}^{\prime
\prime }\right\} \\
&&+\mu ^{\prime }\left( l\right) \left\{ \left( \frac{2f\left( s\right) }{%
\left( \varepsilon _{2}\kappa _{g}^{2}+\varepsilon _{3}\kappa
_{n}^{2}\right) ^{2}}\left[ \left( \varepsilon _{1}\varepsilon _{2}\kappa
_{g}^{2}+\varepsilon _{1}\varepsilon _{3}\kappa _{n}^{2}\right) \left(
-\varepsilon _{3}\kappa _{n}^{\prime }-4\varepsilon _{2}\varepsilon
_{3}\kappa _{g}\tau _{g}\right) \right. \right. \right. \\
&&\left. \left. -2\varepsilon _{1}\kappa _{g}\left( -\varepsilon _{2}\kappa
_{g}\left( \varepsilon _{3}\kappa _{n}^{\prime }+\varepsilon _{2}\varepsilon
_{3}\kappa _{g}\tau _{g}\right) +\varepsilon _{3}\kappa _{n}\left(
\varepsilon _{2}\kappa _{g}^{\prime }-\varepsilon _{2}\varepsilon _{3}\kappa
_{n}\tau _{g}\right) \right) \right] \right) _{s=l} \\
&&\left. -\left( \frac{2\varepsilon _{3}\kappa _{n}f\left( s\right) }{\left(
\varepsilon _{2}\kappa _{g}^{2}+\varepsilon _{3}\kappa _{n}^{2}\right) ^{2}}%
\right) _{s=l}^{\prime }\right\} +\mu ^{\prime \prime }\left( l\right)
\left( \frac{2f\left( s\right) }{\left( \varepsilon _{2}\kappa
_{g}^{2}+\varepsilon _{3}\kappa _{n}^{2}\right) ^{2}}\right) _{s=l} \\
&&-\mu ^{\prime \prime }\left( 0\right) \left( \frac{2f\left( s\right) }{%
\left( \varepsilon _{2}\kappa _{g}^{2}+\varepsilon _{3}\kappa
_{n}^{2}\right) ^{2}}\right) _{s=0}.
\end{eqnarray*}

\subsection{Intrinsic equations for a generalized elastic line on a
spacelike surface}

On a spacelike surface $n$ is timelike. So, we have $T$ and $Q$ are
spacelike and $\varepsilon _{1}=\left\langle T,T\right\rangle =1$, $%
\varepsilon _{2}=\left\langle Q,Q\right\rangle =1$,$\ \varepsilon _{3}=$ $%
\left\langle n,n\right\rangle =-1.\ $Hence $H^{\prime }\left( 0\right) \ $%
can be written as%
\begin{eqnarray*}
H^{\prime }\left( 0\right) &=&\int_{0}^{l}\mu \left\{ \kappa _{g}f^{2}\left(
l\right) +\frac{2f\left( s\right) }{\left( \kappa _{g}^{2}-\kappa
_{n}^{2}\right) ^{2}}\left( \left( \kappa _{g}^{2}-\kappa _{n}^{2}\right) %
\left[ \left( -3\kappa _{g}^{2}\kappa _{n}^{\prime }\right. \right. \right.
\right. \\
&&-3\kappa _{g}^{3}\tau _{g}+3\kappa _{g}\kappa _{g}^{\prime }\kappa
_{n}+2\kappa _{g}\kappa _{n}^{2}\tau _{g}+\kappa _{g}^{4}-\kappa
_{g}^{2}\kappa _{n}^{2}-\kappa _{n}^{\prime }\tau _{g}^{2}-\kappa
_{g}^{\prime }\tau _{g}^{\prime } \\
&&\left. +4\kappa _{n}\tau _{g}\tau _{g}^{\prime }+\kappa _{g}\tau
_{g}^{\prime \prime }+\kappa _{g}\kappa _{n}^{2}\tau _{g}\right) -\left(
\kappa _{g}\left( \kappa _{n}^{\prime }+\kappa _{g}\tau _{g}\right) \right.
\\
&&\left. \left. \left. -\kappa _{n}\left( \kappa _{g}^{\prime }+\kappa
_{n}\tau _{g}\right) \right) \left( -4\kappa _{g}^{3}+4\kappa _{g}\kappa
_{n}^{2}+2\kappa _{g}\tau _{g}^{2}-\kappa _{n}\tau _{g}^{\prime }\right) 
\right] \right) \\
&&-\left[ \frac{2f\left( s\right) }{\left( \kappa _{g}^{2}-\kappa
_{n}^{2}\right) ^{2}}\left[ \left( \kappa _{g}^{2}-\kappa _{n}^{2}\right)
\left( -4\kappa _{g}^{2}\kappa _{n}-\kappa _{n}^{3}-2\kappa _{g}^{\prime
}\tau _{g}\right. \right. \right. \\
&&\left. +2\kappa _{n}\tau _{g}^{2}+3\kappa _{g}\tau _{g}^{\prime }+3\kappa
_{g}^{2}\kappa _{n}+3\kappa _{n}\tau _{g}^{2}\right) -4\kappa _{n}\tau
_{g}\left( \kappa _{g}\left( \kappa _{n}^{\prime }+\kappa _{g}\tau
_{g}\right) \right. \\
&&\left. \left. \left. -\kappa _{n}\left( \kappa _{g}^{\prime }+\kappa
_{n}\tau _{g}\right) \right) \right] \right] ^{\prime }+\left[ \frac{%
2f\left( s\right) }{\left( \kappa _{g}^{2}-\kappa _{n}^{2}\right) ^{2}}\left[
\left( \kappa _{g}^{2}-\kappa _{n}^{2}\right) \left( \kappa _{n}^{\prime
}\right. \right. \right. \\
&&\left. \left. \left. +4\kappa _{g}\tau _{g}\right) -2\kappa _{g}\left(
\kappa _{g}\left( \kappa _{n}^{\prime }+\kappa _{g}\tau _{g}\right) -\kappa
_{n}\left( \kappa _{g}^{\prime }+\kappa _{n}\tau _{g}\right) \right) \right] %
\right] ^{\prime \prime } \\
&&\left. +\left[ \frac{2\kappa _{n}f\left( s\right) }{\left( \kappa
_{g}^{2}-\kappa _{n}^{2}\right) ^{2}}\right] ^{\prime \prime \prime
}\right\} ds \\
&&+\mu \left( l\right) \left\{ \frac{2f\left( l\right) }{\left( \kappa
_{g}^{2}\left( l\right) -\kappa _{n}^{2}\left( l\right) \right) ^{2}}\left[
\left( \kappa _{g}^{2}-\kappa _{n}^{2}\right) \left( -4\kappa _{g}^{2}\kappa
_{n}-\kappa _{n}^{3}-2\kappa _{g}^{\prime }\tau _{g}\right. \right. \right.
\\
&&\left. +2\kappa _{n}\tau _{g}^{2}+3\kappa _{g}\tau _{g}^{\prime }+3\kappa
_{g}^{2}\kappa _{n}+3\kappa _{n}\tau _{g}^{2}\right) -4\kappa _{n}\tau
_{g}\left( \kappa _{g}\left( \kappa _{n}^{\prime }+\kappa _{g}\tau
_{g}\right) \right. \\
&&\left. \left. -\kappa _{n}\left( \kappa _{g}^{\prime }+\kappa _{n}\tau
_{g}\right) \right) \right] _{s=l}-\left( \frac{2f\left( s\right) }{\left(
\kappa _{g}^{2}-\kappa _{n}^{2}\right) ^{2}}\left[ \left( \kappa
_{g}^{2}-\kappa _{n}^{2}\right) \left( \kappa _{n}^{\prime }+4\kappa
_{g}\tau _{g}\right) \right. \right. \\
&&\left. \left. \left. -2\kappa _{g}\left( \kappa _{g}\left( \kappa
_{n}^{\prime }+\kappa _{g}\tau _{g}\right) -\kappa _{n}\left( \kappa
_{g}^{\prime }+\kappa _{n}\tau _{g}\right) \right) \right] \right)
_{s=l}^{\prime }-\left( \frac{2\kappa _{n}f\left( s\right) }{\left( \kappa
_{g}^{2}-\kappa _{n}^{2}\right) ^{2}}\right) _{s=l}^{\prime \prime }\right\}
\\
&&+\mu ^{\prime }\left( l\right) \left\{ \left( \frac{2f\left( s\right) }{%
\left( \kappa _{g}^{2}-\kappa _{n}^{2}\right) ^{2}}\left[ \left( \kappa
_{g}^{2}-\kappa _{n}^{2}\right) \left( \kappa _{n}^{\prime }+4\kappa
_{g}\tau _{g}\right) \right. \right. \right. \\
&&\left. \left. -2\kappa _{g}\left( \kappa _{g}\left( \kappa _{n}^{\prime
}+\kappa _{g}\tau _{g}\right) -\kappa _{n}\left( \kappa _{g}^{\prime
}+\kappa _{n}\tau _{g}\right) \right) \right] \right) _{s=l} \\
&&\left. +\left( \frac{2\kappa _{n}f\left( s\right) }{\left( \kappa
_{g}^{2}-\kappa _{n}^{2}\right) ^{2}}\right) _{s=l}^{\prime }\right\} +\mu
^{\prime \prime }\left( l\right) \left( \frac{2f\left( s\right) }{\left(
\kappa _{g}^{2}-\kappa _{n}^{2}\right) ^{2}}\right) _{s=l} \\
&&-\mu ^{\prime \prime }\left( 0\right) \left( \frac{2f\left( s\right) }{%
\left( \kappa _{g}^{2}-\kappa _{n}^{2}\right) ^{2}}\right) _{s=0}.
\end{eqnarray*}%
In order that $H^{\prime }(0)=0$ for all choices of the function $\mu (s)$
satisfying $\left( \ref{2.5}\right) $, with arbitrary values of $\mu (l)$,$\
\mu ^{\prime }(l)$ and $\mu ^{\prime \prime }(l)$, spacelike arc $\alpha $
must satisfy four boundary conditions%
\begin{eqnarray*}
&&\frac{2f\left( l\right) }{\left( \kappa _{g}^{2}\left( l\right) -\kappa
_{n}^{2}\left( l\right) \right) ^{2}}\left[ \left( \kappa _{g}^{2}-\kappa
_{n}^{2}\right) \left( -4\kappa _{g}^{2}\kappa _{n}-\kappa _{n}^{3}-2\kappa
_{g}^{\prime }\tau _{g}\right. \right. \\
&&\left. +2\kappa _{n}\tau _{g}^{2}+3\kappa _{g}\tau _{g}^{\prime }+3\kappa
_{g}^{2}\kappa _{n}+3\kappa _{n}\tau _{g}^{2}\right) -4\kappa _{n}\tau
_{g}\left( \kappa _{g}\left( \kappa _{n}^{\prime }+\kappa _{g}\tau
_{g}\right) \right. \\
&&\left. \left. -\kappa _{n}\left( \kappa _{g}^{\prime }+\kappa _{n}\tau
_{g}\right) \right) \right] _{s=l}-\left( \frac{2f\left( s\right) }{\left(
\kappa _{g}^{2}-\kappa _{n}^{2}\right) ^{2}}\left[ \left( \kappa
_{g}^{2}-\kappa _{n}^{2}\right) \left( \kappa _{n}^{\prime }+4\kappa
_{g}\tau _{g}\right) \right. \right. \\
&&\left. \left. -2\kappa _{g}\left( \kappa _{g}\left( \kappa _{n}^{\prime
}+\kappa _{g}\tau _{g}\right) -\kappa _{n}\left( \kappa _{g}^{\prime
}+\kappa _{n}\tau _{g}\right) \right) \right] \right) _{s=l}^{\prime
}-\left( \frac{2\kappa _{n}f\left( s\right) }{\left( \kappa _{g}^{2}-\kappa
_{n}^{2}\right) ^{2}}\right) _{s=l}^{\prime \prime }=0,
\end{eqnarray*}%
\begin{eqnarray*}
&&\left( \frac{2f\left( s\right) }{\left( \kappa _{g}^{2}-\kappa
_{n}^{2}\right) ^{2}}\left[ \left( \kappa _{g}^{2}-\kappa _{n}^{2}\right)
\left( \kappa _{n}^{\prime }+4\kappa _{g}\tau _{g}\right) \right. \right. \\
&&\left. \left. -2\kappa _{g}\left( \kappa _{g}\left( \kappa _{n}^{\prime
}+\kappa _{g}\tau _{g}\right) -\kappa _{n}\left( \kappa _{g}^{\prime
}+\kappa _{n}\tau _{g}\right) \right) \right] \right) _{s=l} \\
&&+\left( \frac{2\kappa _{n}f\left( s\right) }{\left( \kappa _{g}^{2}-\kappa
_{n}^{2}\right) ^{2}}\right) _{s=l}^{\prime }=0,\ \ \ \ \ \ \ \ \ \ \ \ \ \
\ \ \ \ \ \ \ \ \ \ \ \ \ \ \ \ \ \ \ \ \ \ \ \ \ \ \ \ \ \ \ \ \ \ \ \ \ \
\ \ \ \ \ \ \ \ \ \ \ \ \ \ \ \ \ \ \ \ \ 
\end{eqnarray*}%
\begin{equation*}
f\left( l\right) =0,\ \ \ \ \ \ \ \ \ \ \ \ \ \ \ \ \ \ \ \ \ \ \ \ \ \ \ \
\ \ \ \ \ \ \ \ \ \ \ \ \ \ \ \ \ \ \ \ \ \ \ \ \ \ \ \ \ \ \ \ \ \ \ \ \ \
\ \ \ \ \ \ \ \ \ \ \ 
\end{equation*}%
\begin{equation*}
f\left( 0\right) =0,\ \ \ \ \ \ \ \ \ \ \ \ \ \ \ \ \ \ \ \ \ \ \ \ \ \ \ \
\ \ \ \ \ \ \ \ \ \ \ \ \ \ \ \ \ \ \ \ \ \ \ \ \ \ \ \ \ \ \ \ \ \ \ \ \ \
\ \ \ \ \ \ \ \ \ \ 
\end{equation*}%
and differential equation%
\begin{eqnarray*}
&&\kappa _{g}f^{2}\left( l\right) +\frac{2f\left( s\right) }{\left( \kappa
_{g}^{2}-\kappa _{n}^{2}\right) ^{2}}\left( \left( \kappa _{g}^{2}-\kappa
_{n}^{2}\right) \left[ \left( -3\kappa _{g}^{2}\kappa _{n}^{\prime }-3\kappa
_{g}^{3}\tau _{g}\right. \right. \right. \\
&&+3\kappa _{g}\kappa _{g}^{\prime }\kappa _{n}+2\kappa _{g}\kappa
_{n}^{2}\tau _{g}+\kappa _{g}^{4}-\kappa _{g}^{2}\kappa _{n}^{2}-\kappa
_{n}^{\prime }\tau _{g}^{2}-\kappa _{g}^{\prime }\tau _{g}^{\prime } \\
&&\left. +4\kappa _{n}\tau _{g}\tau _{g}^{\prime }+\kappa _{g}\tau
_{g}^{\prime \prime }+\kappa _{g}\kappa _{n}^{2}\tau _{g}\right) -\left(
\kappa _{g}\left( \kappa _{n}^{\prime }+\kappa _{g}\tau _{g}\right) \right.
\\
&&\left. \left. \left. -\kappa _{n}\left( \kappa _{g}^{\prime }+\kappa
_{n}\tau _{g}\right) \right) \left( -4\kappa _{g}^{3}+4\kappa _{g}\kappa
_{n}^{2}+2\kappa _{g}\tau _{g}^{2}-\kappa _{n}\tau _{g}^{\prime }\right) 
\right] \right) \\
&&-\left[ \frac{2f\left( s\right) }{\left( \kappa _{g}^{2}-\kappa
_{n}^{2}\right) ^{2}}\left[ \left( \kappa _{g}^{2}-\kappa _{n}^{2}\right)
\left( -4\kappa _{g}^{2}\kappa _{n}-\kappa _{n}^{3}-2\kappa _{g}^{\prime
}\tau _{g}\right. \right. \right. \\
&&\left. +2\kappa _{n}\tau _{g}^{2}+3\kappa _{g}\tau _{g}^{\prime }+3\kappa
_{g}^{2}\kappa _{n}+3\kappa _{n}\tau _{g}^{2}\right) -4\kappa _{n}\tau
_{g}\left( \kappa _{g}\left( \kappa _{n}^{\prime }+\kappa _{g}\tau
_{g}\right) \right. \\
&&\left. \left. \left. -\kappa _{n}\left( \kappa _{g}^{\prime }+\kappa
_{n}\tau _{g}\right) \right) \right] \right] ^{\prime }+\left[ \frac{%
2f\left( s\right) }{\left( \kappa _{g}^{2}-\kappa _{n}^{2}\right) ^{2}}\left[
\left( \kappa _{g}^{2}-\kappa _{n}^{2}\right) \left( \kappa _{n}^{\prime
}\right. \right. \right. \\
&&\left. \left. \left. +4\kappa _{g}\tau _{g}\right) -2\kappa _{g}\left(
\kappa _{g}\left( \kappa _{n}^{\prime }+\kappa _{g}\tau _{g}\right) -\kappa
_{n}\left( \kappa _{g}^{\prime }+\kappa _{n}\tau _{g}\right) \right) \right] %
\right] ^{\prime \prime } \\
&&+\left[ \frac{2\kappa _{n}f\left( s\right) }{\left( \kappa _{g}^{2}-\kappa
_{n}^{2}\right) ^{2}}\right] ^{\prime \prime \prime }=0,
\end{eqnarray*}%
where%
\begin{equation*}
f\left( s\right) =\frac{\kappa _{g}\left( \kappa _{n}^{\prime }+\kappa
_{g}\tau _{g}\right) -\kappa _{n}\left( \kappa _{g}^{\prime }+\kappa
_{n}\tau _{g}\right) }{\kappa _{g}^{2}-\kappa _{n}^{2}}.
\end{equation*}

\subsection{Intrinsic equations for a generalized elastic line on a timelike
surface for timelike arc $\protect\alpha $}

Since $\alpha \ $is timelike $T$ is timelike. So, we have $Q$ and $n$ are
spacelike and $\varepsilon _{1}=\left\langle T,T\right\rangle =-1$, $%
\varepsilon _{2}=\left\langle Q,Q\right\rangle =1$,$\ \varepsilon _{3}=$ $%
\left\langle n,n\right\rangle =1.\ $Hence $H^{\prime }\left( 0\right) \ $can
be written as%
\begin{eqnarray*}
H^{\prime }\left( 0\right) &=&\int_{0}^{l}\mu \left\{ -\kappa
_{g}f^{2}\left( l\right) +\frac{2f\left( s\right) }{\left( \kappa
_{g}^{2}+\kappa _{n}^{2}\right) ^{2}}\left( \left( -\kappa _{g}^{2}-\kappa
_{n}^{2}\right) \left[ \left( -3\kappa _{g}^{2}\kappa _{n}^{\prime }\right.
\right. \right. \right. \\
&&-3\kappa _{g}^{3}\tau _{g}+3\kappa _{g}\kappa _{g}^{\prime }\kappa
_{n}-2\kappa _{g}\kappa _{n}^{2}\tau _{g}+\kappa _{g}^{4}+\kappa
_{g}^{2}\kappa _{n}^{2}-\kappa _{n}^{\prime }\tau _{g}^{2}+\kappa
_{g}^{\prime }\tau _{g}^{\prime } \\
&&\left. +4\kappa _{n}\tau _{g}\tau _{g}^{\prime }-\kappa _{g}\tau
_{g}^{\prime \prime }-\kappa _{g}\kappa _{n}^{2}\tau _{g}\right) -\left(
-\kappa _{g}\left( \kappa _{n}^{\prime }+\kappa _{g}\tau _{g}\right) \right.
\\
&&\left. \left. \left. +\kappa _{n}\left( \kappa _{g}^{\prime }-\kappa
_{n}\tau _{g}\right) \right) \left( -4\kappa _{g}^{3}-4\kappa _{g}\kappa
_{n}^{2}+2\kappa _{g}\tau _{g}^{2}-2\kappa _{n}\tau _{g}^{\prime }\right) 
\right] \right) \\
&&-\left[ \frac{2f\left( s\right) }{\left( \kappa _{g}^{2}+\kappa
_{n}^{2}\right) ^{2}}\left[ \left( -\kappa _{g}^{2}-\kappa _{n}^{2}\right)
\left( -4\kappa _{g}^{2}\kappa _{n}+\kappa _{n}^{3}+2\kappa _{g}^{\prime
}\tau _{g}\right. \right. \right. \\
&&\left. +2\kappa _{n}\tau _{g}^{2}-3\kappa _{g}\tau _{g}^{\prime }+3\kappa
_{g}^{2}\kappa _{n}+3\kappa _{n}\tau _{g}^{2}\right) +4\kappa _{n}\tau
_{g}\left( -\kappa _{g}\left( \kappa _{n}^{\prime }+\kappa _{g}\tau
_{g}\right) \right. \\
&&\left. \left. \left. +\kappa _{n}\left( \kappa _{g}^{\prime }-\kappa
_{n}\tau _{g}\right) \right) \right] \right] ^{\prime }+\left[ \frac{%
2f\left( s\right) }{\left( \kappa _{g}^{2}+\kappa _{n}^{2}\right) ^{2}}\left[
\left( -\kappa _{g}^{2}-\kappa _{n}^{2}\right) \left( -\kappa _{n}^{\prime
}\right. \right. \right. \\
&&\left. \left. \left. -4\kappa _{g}\tau _{g}\right) +2\kappa _{g}\left(
-\kappa _{g}\left( \kappa _{n}^{\prime }+\kappa _{g}\tau _{g}\right) +\kappa
_{n}\left( \kappa _{g}^{\prime }-\kappa _{n}\tau _{g}\right) \right) \right] %
\right] ^{\prime \prime } \\
&&\left. -\left[ \frac{2\kappa _{n}f\left( s\right) }{\left( \kappa
_{g}^{2}+\kappa _{n}^{2}\right) ^{2}}\right] ^{\prime \prime \prime
}\right\} ds \\
&&+\mu \left( l\right) \left\{ \frac{2f\left( l\right) }{\left( \kappa
_{g}^{2}\left( l\right) +\kappa _{n}^{2}\left( l\right) \right) ^{2}}\left[
\left( -\kappa _{g}^{2}-\kappa _{n}^{2}\right) \left( -4\kappa
_{g}^{2}\kappa _{n}+\kappa _{n}^{3}+2\kappa _{g}^{\prime }\tau _{g}\right.
\right. \right. \\
&&\left. +2\kappa _{n}\tau _{g}^{2}-3\kappa _{g}\tau _{g}^{\prime }+3\kappa
_{g}^{2}\kappa _{n}+3\kappa _{n}\tau _{g}^{2}\right) +4\kappa _{n}\tau
_{g}\left( -\kappa _{g}\left( \kappa _{n}^{\prime }+\kappa _{g}\tau
_{g}\right) \right. \\
&&\left. \left. +\kappa _{n}\left( \kappa _{g}^{\prime }-\kappa _{n}\tau
_{g}\right) \right) \right] _{s=l}-\left( \frac{2f\left( s\right) }{\left(
\kappa _{g}^{2}+\kappa _{n}^{2}\right) ^{2}}\left[ \left( -\kappa
_{g}^{2}-\kappa _{n}^{2}\right) \left( -\kappa _{n}^{\prime }-4\kappa
_{g}\tau _{g}\right) \right. \right. \\
&&\left. \left. \left. \left. +2\kappa _{g}\left( -\kappa _{g}\left( \kappa
_{n}^{\prime }+\kappa _{g}\tau _{g}\right) +\kappa _{n}\left( \kappa
_{g}^{\prime }-\kappa _{n}\tau _{g}\right) \right) \right) \right] \right)
_{s=l}^{\prime }+\left( \frac{2\kappa _{n}f\left( s\right) }{\left( \kappa
_{g}^{2}+\kappa _{n}^{2}\right) ^{2}}\right) _{s=l}^{\prime \prime }\right\}
\\
&&+\mu ^{\prime }\left( l\right) \left\{ \left( \frac{2f\left( s\right) }{%
\left( \kappa _{g}^{2}+\kappa _{n}^{2}\right) ^{2}}\left[ \left( -\kappa
_{g}^{2}-\kappa _{n}^{2}\right) \left( -\kappa _{n}^{\prime }-4\kappa
_{g}\tau _{g}\right) \right. \right. \right. \\
&&\left. \left. +\kappa _{g}\left( -\kappa _{g}\left( \kappa _{n}^{\prime
}+\kappa _{g}\tau _{g}\right) +\kappa _{n}\left( \kappa _{g}^{\prime
}-\kappa _{n}\tau _{g}\right) \right) \right] \right) _{s=l} \\
&&\left. -\left( \frac{2\kappa _{n}f\left( s\right) }{\left( \kappa
_{g}^{2}+\kappa _{n}^{2}\right) ^{2}}\right) _{s=l}^{\prime }\right\} +\mu
^{\prime \prime }\left( l\right) \left( \frac{2f\left( s\right) }{\left(
\kappa _{g}^{2}+\kappa _{n}^{2}\right) ^{2}}\right) _{s=l} \\
&&-\mu ^{\prime \prime }\left( 0\right) \left( \frac{2f\left( s\right) }{%
\left( \kappa _{g}^{2}+\kappa _{n}^{2}\right) ^{2}}\right) _{s=0}.
\end{eqnarray*}%
In order that $H^{\prime }(0)=0$ for all choices of the function $\mu (s)$
satisfying $\left( \ref{2.5}\right) $, with arbitrary values of $\mu (l)$,$\
\mu ^{\prime }(l)$ and $\mu ^{\prime \prime }(l)$, spacelike arc $\alpha $
must satisfy four boundary conditions%
\begin{eqnarray*}
&&\frac{2f\left( l\right) }{\left( \kappa _{g}^{2}\left( l\right) +\kappa
_{n}^{2}\left( l\right) \right) ^{2}}\left[ \left( -\kappa _{g}^{2}-\kappa
_{n}^{2}\right) \left( -4\kappa _{g}^{2}\kappa _{n}+\kappa _{n}^{3}+2\kappa
_{g}^{\prime }\tau _{g}\right. \right. \\
&&\left. +2\kappa _{n}\tau _{g}^{2}-3\kappa _{g}\tau _{g}^{\prime }+3\kappa
_{g}^{2}\kappa _{n}+3\kappa _{n}\tau _{g}^{2}\right) +4\kappa _{n}\tau
_{g}\left( -\kappa _{g}\left( \kappa _{n}^{\prime }+\kappa _{g}\tau
_{g}\right) \right. \\
&&\left. \left. +\kappa _{n}\left( \kappa _{g}^{\prime }-\kappa _{n}\tau
_{g}\right) \right) \right] _{s=l}-\left( \frac{2f\left( s\right) }{\left(
\kappa _{g}^{2}+\kappa _{n}^{2}\right) ^{2}}\left[ \left( -\kappa
_{g}^{2}-\kappa _{n}^{2}\right) \left( -\kappa _{n}^{\prime }-4\kappa
_{g}\tau _{g}\right) \right. \right. \\
&&\left. \left. \left. +2\kappa _{g}\left( -\kappa _{g}\left( \kappa
_{n}^{\prime }+\kappa _{g}\tau _{g}\right) +\kappa _{n}\left( \kappa
_{g}^{\prime }-\kappa _{n}\tau _{g}\right) \right) \right) \right] \right)
_{s=l}^{\prime }+\left( \frac{2\kappa _{n}f\left( s\right) }{\left( \kappa
_{g}^{2}+\kappa _{n}^{2}\right) ^{2}}\right) _{s=l}^{\prime \prime }=0,\ \ \
\ \ \ \ \ \ \ \ \ \ \ \ \ \ \ \ 
\end{eqnarray*}%
\begin{eqnarray*}
&&\left( \frac{2f\left( s\right) }{\left( \kappa _{g}^{2}+\kappa
_{n}^{2}\right) ^{2}}\left[ \left( -\kappa _{g}^{2}-\kappa _{n}^{2}\right)
\left( -\kappa _{n}^{\prime }-4\kappa _{g}\tau _{g}\right) +\kappa
_{g}\left( -\kappa _{g}\left( \kappa _{n}^{\prime }+\kappa _{g}\tau
_{g}\right) +\kappa _{n}\left( \kappa _{g}^{\prime }-\kappa _{n}\tau
_{g}\right) \right) \right] \right) _{s=l} \\
&&-\left( \frac{2\kappa _{n}f\left( s\right) }{\left( \kappa _{g}^{2}+\kappa
_{n}^{2}\right) ^{2}}\right) _{s=l}^{\prime }=0,
\end{eqnarray*}%
\begin{equation*}
f\left( l\right) =0,\ \ \ \ \ \ \ \ \ \ \ \ \ \ \ \ \ \ \ \ \ \ \ \ \ \ \ \
\ \ \ \ \ \ \ \ \ \ \ \ \ \ \ \ \ \ \ \ \ \ \ \ \ \ \ \ \ \ \ \ \ \ \ \ \ \
\ \ \ \ \ \ \ \ \ \ \ 
\end{equation*}%
\begin{equation*}
f\left( 0\right) =0,\ \ \ \ \ \ \ \ \ \ \ \ \ \ \ \ \ \ \ \ \ \ \ \ \ \ \ \
\ \ \ \ \ \ \ \ \ \ \ \ \ \ \ \ \ \ \ \ \ \ \ \ \ \ \ \ \ \ \ \ \ \ \ \ \ \
\ \ \ \ \ \ \ \ \ \ 
\end{equation*}%
and differential equation%
\begin{eqnarray*}
&&-\kappa _{g}f^{2}\left( l\right) +\frac{2f\left( s\right) }{\left( \kappa
_{g}^{2}+\kappa _{n}^{2}\right) ^{2}}\left( \left( -\kappa _{g}^{2}-\kappa
_{n}^{2}\right) \left[ \left( -3\kappa _{g}^{2}\kappa _{n}^{\prime }\right.
\right. \right. \\
&&-3\kappa _{g}^{3}\tau _{g}+3\kappa _{g}\kappa _{g}^{\prime }\kappa
_{n}-2\kappa _{g}\kappa _{n}^{2}\tau _{g}+\kappa _{g}^{4}+\kappa
_{g}^{2}\kappa _{n}^{2}-\kappa _{n}^{\prime }\tau _{g}^{2}+\kappa
_{g}^{\prime }\tau _{g}^{\prime } \\
&&\left. +4\kappa _{n}\tau _{g}\tau _{g}^{\prime }-\kappa _{g}\tau
_{g}^{\prime \prime }-\kappa _{g}\kappa _{n}^{2}\tau _{g}\right) -\left(
-\kappa _{g}\left( \kappa _{n}^{\prime }+\kappa _{g}\tau _{g}\right) \right.
\\
&&\left. \left. \left. +\kappa _{n}\left( \kappa _{g}^{\prime }-\kappa
_{n}\tau _{g}\right) \right) \left( -4\kappa _{g}^{3}-4\kappa _{g}\kappa
_{n}^{2}+2\kappa _{g}\tau _{g}^{2}-2\kappa _{n}\tau _{g}^{\prime }\right) 
\right] \right) \\
&&-\left[ \frac{2f\left( s\right) }{\left( \kappa _{g}^{2}+\kappa
_{n}^{2}\right) ^{2}}\left[ \left( -\kappa _{g}^{2}-\kappa _{n}^{2}\right)
\left( -4\kappa _{g}^{2}\kappa _{n}+\kappa _{n}^{3}+2\kappa _{g}^{\prime
}\tau _{g}\right. \right. \right. \\
&&\left. +2\kappa _{n}\tau _{g}^{2}-3\kappa _{g}\tau _{g}^{\prime }+3\kappa
_{g}^{2}\kappa _{n}+3\kappa _{n}\tau _{g}^{2}\right) +4\kappa _{n}\tau
_{g}\left( -\kappa _{g}\left( \kappa _{n}^{\prime }+\kappa _{g}\tau
_{g}\right) \right. \\
&&\left. \left. \left. +\kappa _{n}\left( \kappa _{g}^{\prime }-\kappa
_{n}\tau _{g}\right) \right) \right] \right] ^{\prime }+\left[ \frac{%
2f\left( s\right) }{\left( \kappa _{g}^{2}+\kappa _{n}^{2}\right) ^{2}}\left[
\left( -\kappa _{g}^{2}-\kappa _{n}^{2}\right) \left( -\kappa _{n}^{\prime
}\right. \right. \right. \\
&&\left. \left. \left. -4\kappa _{g}\tau _{g}\right) +2\kappa _{g}\left(
-\kappa _{g}\left( \kappa _{n}^{\prime }+\kappa _{g}\tau _{g}\right) +\kappa
_{n}\left( \kappa _{g}^{\prime }-\kappa _{n}\tau _{g}\right) \right) \right] %
\right] ^{\prime \prime } \\
&&-\left[ \frac{2\kappa _{n}f\left( s\right) }{\left( \kappa _{g}^{2}+\kappa
_{n}^{2}\right) ^{2}}\right] ^{\prime \prime \prime }=0,
\end{eqnarray*}%
where%
\begin{equation*}
f\left( s\right) =\frac{\kappa _{g}\left( \kappa _{n}^{\prime }+\kappa
_{g}\tau _{g}\right) -\kappa _{n}\left( \kappa _{g}^{\prime }-\kappa
_{n}\tau _{g}\right) }{\kappa _{g}^{2}+\kappa _{n}^{2}}.
\end{equation*}

\subsection{Intrinsic equations for a generalized elastic line on a timelike
surface for spacelike arc $\protect\alpha $}

Now $Q$ is timelike, $T$ and $n$ are spacelike. So we have $\varepsilon
_{1}=\left\langle T,T\right\rangle =1$, $\varepsilon _{2}=\left\langle
Q,Q\right\rangle =-1$,$\ \varepsilon _{3}=$ $\left\langle n,n\right\rangle
=1.\ $Hence $H^{\prime }\left( 0\right) \ $can be written as%
\begin{eqnarray*}
H^{\prime }\left( 0\right) &=&\int_{0}^{l}\mu \left\{ \kappa _{g}f^{2}\left(
l\right) +\frac{2f\left( s\right) }{\left( -\kappa _{g}^{2}+\kappa
_{n}^{2}\right) ^{2}}\left( \left( -\kappa _{g}^{2}+\kappa _{n}^{2}\right) %
\left[ \left( -3\kappa _{g}^{2}\kappa _{n}^{\prime }\right. \right. \right.
\right. \\
&&+3\kappa _{g}^{3}\tau _{g}+3\kappa _{g}\kappa _{g}^{\prime }\kappa
_{n}-2\kappa _{g}\kappa _{n}^{2}\tau _{g}-\kappa _{g}^{4}+\kappa
_{g}^{2}\kappa _{n}^{2}+\kappa _{n}^{\prime }\tau _{g}^{2}-\kappa
_{g}^{\prime }\tau _{g}^{\prime } \\
&&\left. +4\kappa _{n}\tau _{g}\tau _{g}^{\prime }+\kappa _{g}\tau
_{g}^{\prime \prime }-\kappa _{g}\kappa _{n}^{2}\tau _{g}\right) -\left(
\kappa _{g}\left( \kappa _{n}^{\prime }-\kappa _{g}\tau _{g}\right) \right.
\\
&&\left. \left. \left. +\kappa _{n}\left( -\kappa _{g}^{\prime }+\kappa
_{n}\tau _{g}\right) \right) \left( 4\kappa _{g}^{3}-4\kappa _{g}\kappa
_{n}^{2}+2\kappa _{g}\tau _{g}^{2}+2\kappa _{n}\tau _{g}^{\prime }\right) 
\right] \right) \\
&&-\left[ \frac{2f\left( s\right) }{\left( -\kappa _{g}^{2}+\kappa
_{n}^{2}\right) ^{2}}\left[ \left( -\kappa _{g}^{2}+\kappa _{n}^{2}\right)
\left( -4\kappa _{g}^{2}\kappa _{n}-\kappa _{n}^{3}-2\kappa _{g}^{\prime
}\tau _{g}\right. \right. \right. \\
&&\left. -2\kappa _{n}\tau _{g}^{2}+3\kappa _{g}\tau _{g}^{\prime }+3\kappa
_{g}^{2}\kappa _{n}-3\kappa _{n}\tau _{g}^{2}\right) +4\kappa _{n}\tau
_{g}\left( \kappa _{g}\left( \kappa _{n}^{\prime }-\kappa _{g}\tau
_{g}\right) \right. \\
&&\left. \left. \left. +\kappa _{n}\left( -\kappa _{g}^{\prime }+\kappa
_{n}\tau _{g}\right) \right) \right] \right] ^{\prime }+\left[ \frac{%
2f\left( s\right) }{\left( -\kappa _{g}^{2}+\kappa _{n}^{2}\right) ^{2}}%
\left[ \left( -\kappa _{g}^{2}+\kappa _{n}^{2}\right) \left( -\kappa
_{n}^{\prime }\right. \right. \right. \\
&&\left. \left. \left. +4\kappa _{g}\tau _{g}\right) -2\kappa _{g}\left(
\kappa _{g}\left( \kappa _{n}^{\prime }-\kappa _{g}\tau _{g}\right) +\kappa
_{n}\left( -\kappa _{g}^{\prime }+\kappa _{n}\tau _{g}\right) \right) \right]
\right] ^{\prime \prime } \\
&&\left. -\left[ \frac{2\kappa _{n}f\left( s\right) }{\left( -\kappa
_{g}^{2}+\kappa _{n}^{2}\right) ^{2}}\right] ^{\prime \prime \prime
}\right\} ds \\
&&+\mu \left( l\right) \left\{ \frac{2f\left( l\right) }{\left( -\kappa
_{g}^{2}\left( l\right) +\kappa _{n}^{2}\left( l\right) \right) ^{2}}\left[
\left( -\kappa _{g}^{2}+\kappa _{n}^{2}\right) \left( -4\kappa
_{g}^{2}\kappa _{n}-\kappa _{n}^{3}-2\kappa _{g}^{\prime }\tau _{g}\right.
\right. \right. \\
&&\left. -2\kappa _{n}\tau _{g}^{2}+3\kappa _{g}\tau _{g}^{\prime }+3\kappa
_{g}^{2}\kappa _{n}-3\kappa _{n}\tau _{g}^{2}\right) +4\kappa _{n}\tau
_{g}\left( \kappa _{g}\left( \kappa _{n}^{\prime }-\kappa _{g}\tau
_{g}\right) \right. \\
&&\left. \left. +\kappa _{n}\left( -\kappa _{g}^{\prime }+\kappa _{n}\tau
_{g}\right) \right) \right] _{s=l}-\left( \frac{2f\left( s\right) }{\left(
-\kappa _{g}^{2}+\kappa _{n}^{2}\right) ^{2}}\left[ \left( -\kappa
_{g}^{2}+\kappa _{n}^{2}\right) \left( -\kappa _{n}^{\prime }+4\kappa
_{g}\tau _{g}\right) \right. \right. \\
&&-2\kappa _{g}\left( \kappa _{g}\left( \kappa _{n}^{\prime }-\kappa
_{g}\tau _{g}\right) \right. \\
&&\left. \left. \left. \left. +\kappa _{n}\left( -\kappa _{g}^{\prime
}+\kappa _{n}\tau _{g}\right) \right) \right] \right) _{s=l}^{\prime
}+\left( \frac{2\kappa _{n}f\left( s\right) }{\left( -\kappa _{g}^{2}+\kappa
_{n}^{2}\right) ^{2}}\right) _{s=l}^{\prime \prime }\right\} \\
&&+\mu ^{\prime }\left( l\right) \left\{ \left( \frac{2f\left( s\right) }{%
\left( -\kappa _{g}^{2}+\kappa _{n}^{2}\right) ^{2}}\left[ \left( -\kappa
_{g}^{2}+\kappa _{n}^{2}\right) \left( -\kappa _{n}^{\prime }+4\kappa
_{g}\tau _{g}\right) \right. \right. \right. \\
&&\left. \left. -2\kappa _{g}\left( \kappa _{g}\left( \kappa _{n}^{\prime
}-\kappa _{g}\tau _{g}\right) +\kappa _{n}\left( -\kappa _{g}^{\prime
}+\kappa _{n}\tau _{g}\right) \right) \right] \right) _{s=l} \\
&&\left. -\left( \frac{2\kappa _{n}f\left( s\right) }{\left( -\kappa
_{g}^{2}+\kappa _{n}^{2}\right) ^{2}}\right) _{s=l}^{\prime }\right\} +\mu
^{\prime \prime }\left( l\right) \left( \frac{2f\left( s\right) }{\left(
-\kappa _{g}^{2}+\kappa _{n}^{2}\right) ^{2}}\right) _{s=l} \\
&&-\mu ^{\prime \prime }\left( 0\right) \left( \frac{2f\left( s\right) }{%
\left( -\kappa _{g}^{2}+\kappa _{n}^{2}\right) ^{2}}\right) _{s=0}.
\end{eqnarray*}%
In order that $H^{\prime }(0)=0$ for all choices of the function $\mu (s)$
satisfying $\left( \ref{2.5}\right) $, with arbitrary values of $\mu (l)$,$\
\mu ^{\prime }(l)$ and $\mu ^{\prime \prime }(l)$, spacelike arc $\alpha $
must satisfy four boundary conditions%
\begin{eqnarray*}
&&\frac{2f\left( l\right) }{\left( -\kappa _{g}^{2}\left( l\right) +\kappa
_{n}^{2}\left( l\right) \right) ^{2}}\left[ \left( -\kappa _{g}^{2}+\kappa
_{n}^{2}\right) \left( -4\kappa _{g}^{2}\kappa _{n}-\kappa _{n}^{3}-2\kappa
_{g}^{\prime }\tau _{g}\right. \right. \\
&&\left. -2\kappa _{n}\tau _{g}^{2}+3\kappa _{g}\tau _{g}^{\prime }+3\kappa
_{g}^{2}\kappa _{n}-3\kappa _{n}\tau _{g}^{2}\right) +4\kappa _{n}\tau
_{g}\left( \kappa _{g}\left( \kappa _{n}^{\prime }-\kappa _{g}\tau
_{g}\right) \right. \\
&&\left. \left. +\kappa _{n}\left( -\kappa _{g}^{\prime }+\kappa _{n}\tau
_{g}\right) \right) \right] _{s=l}-\left( \frac{2f\left( s\right) }{\left(
-\kappa _{g}^{2}+\kappa _{n}^{2}\right) ^{2}}\left[ \left( -\kappa
_{g}^{2}+\kappa _{n}^{2}\right) \left( -\kappa _{n}^{\prime }+4\kappa
_{g}\tau _{g}\right) \right. \right. \\
&&-2\kappa _{g}\left( \kappa _{g}\left( \kappa _{n}^{\prime }-\kappa
_{g}\tau _{g}\right) \right. \\
&&\left. \left. \left. +\kappa _{n}\left( -\kappa _{g}^{\prime }+\kappa
_{n}\tau _{g}\right) \right) \right] \right) _{s=l}^{\prime }+\left( \frac{%
2\kappa _{n}f\left( s\right) }{\left( -\kappa _{g}^{2}+\kappa
_{n}^{2}\right) ^{2}}\right) _{s=l}^{\prime \prime }=0,
\end{eqnarray*}%
\begin{eqnarray*}
&&\left( \frac{2f\left( s\right) }{\left( -\kappa _{g}^{2}+\kappa
_{n}^{2}\right) ^{2}}\left[ \left( -\kappa _{g}^{2}+\kappa _{n}^{2}\right)
\left( -\kappa _{n}^{\prime }+4\kappa _{g}\tau _{g}\right) \right. \right. \\
&&\left. \left. -2\kappa _{g}\left( \kappa _{g}\left( \kappa _{n}^{\prime
}-\kappa _{g}\tau _{g}\right) +\kappa _{n}\left( -\kappa _{g}^{\prime
}+\kappa _{n}\tau _{g}\right) \right) \right] \right) _{s=l} \\
&&-\left( \frac{2\kappa _{n}f\left( s\right) }{\left( -\kappa
_{g}^{2}+\kappa _{n}^{2}\right) ^{2}}\right) _{s=l}^{\prime }=0,\ \ \ \ \ \
\ \ \ \ \ \ \ \ \ \ \ \ \ \ \ \ \ \ \ \ \ \ \ \ \ \ \ \ \ \ \ \ \ \ \ \ \ \
\ \ \ \ \ \ \ \ \ \ \ \ \ \ \ \ \ \ \ \ \ \ \ \ 
\end{eqnarray*}%
\begin{equation*}
f\left( l\right) =0,\ \ \ \ \ \ \ \ \ \ \ \ \ \ \ \ \ \ \ \ \ \ \ \ \ \ \ \
\ \ \ \ \ \ \ \ \ \ \ \ \ \ \ \ \ \ \ \ \ \ \ \ \ \ \ \ \ \ \ \ \ \ \ \ \ \
\ \ \ \ \ \ \ \ \ \ \ 
\end{equation*}%
\begin{equation*}
f\left( 0\right) =0,\ \ \ \ \ \ \ \ \ \ \ \ \ \ \ \ \ \ \ \ \ \ \ \ \ \ \ \
\ \ \ \ \ \ \ \ \ \ \ \ \ \ \ \ \ \ \ \ \ \ \ \ \ \ \ \ \ \ \ \ \ \ \ \ \ \
\ \ \ \ \ \ \ \ \ \ 
\end{equation*}%
where%
\begin{equation*}
f\left( s\right) =\frac{\kappa _{g}\left( \kappa _{n}^{\prime }-\kappa
_{g}\tau _{g}\right) -\kappa _{n}\left( \kappa _{g}^{\prime }+\kappa
_{n}\tau _{g}\right) }{-\kappa _{g}^{2}+\kappa _{n}^{2}}.
\end{equation*}

\end{document}